\documentclass{amsart}[12pt]
\usepackage{amssymb}
\usepackage{amsfonts}
\usepackage{latexsym}

\newtheorem{theorem}{Theorem}[section]

\newtheorem{proposition}[theorem]{Proposition}
\newtheorem{corollary}[theorem]{Corollary} 
\theoremstyle{definition}

\newtheorem{conjecture}[theorem]{Conjecture}

\newcommand{\Tr}{\text{Tr}}

\newcommand{\g}{\mathfrak{g}}

\newcommand{\ben}{\begin{enumerate}}
\newcommand{\een}{\end{enumerate}}

\begin{document}

\title{On the Cachazo-Douglas-Seiberg-Witten conjecture for 
simple Lie algebras, II}

\begin{abstract} 
Recently, motivated by supersymmetric gauge theory, 
Cachazo, Douglas, Seiberg, and Witten proposed a conjecture 
about finite dimensional simple Lie
algebras, and checked it in the classical cases [CDSW,W].
Later V. Kac and the author proposed a uniform approach 
to this conjecture [EK], based on the theory of 
abelian ideals in the Borel subalgebra; this allowed them to check 
the conjecture for type $G_2$.

In this note we further develop this approach, 
and propose three natural conjectures which imply the three 
parts of the CDSW conjecture. 
In a sense, these conjectures explain why 
the CDSW conjecture should be true. 
We show that our conjectures hold
for classical Lie algebras and for $G_2$; 
this, in particular, gives a purely algebraic proof 
of the CDSW conjecture for $SO(N)$ 
(the proof in [W] uses the theory of instantons).  
\end{abstract}

\author{Pavel Etingof}
\address{Department of Mathematics, Massachusetts Institute of Technology,
Cambridge, MA 02139, USA}
\email{etingof@math.mit.edu}

\maketitle

\section{The CDSW conjecture}

Let $\g$ be a simple finite dimensional Lie algebra over ${\Bbb C}$.
We fix an invariant form on $\g$ and do not distinguish $\g$ and
$\g^*$. Let $g$ be the dual Coxeter number of $\g$.  
Consider the associative algebra $R=\wedge(\g\oplus \g)$. This algebra 
is naturally $\Bbb Z_+$-bigraded 
(with the two copies of $\g$ sitting in degrees (1,0) and (0,1), 
respectively). The degree (2,0),(1,1), and (0,2) components 
of $R$ are $\wedge^2\g,\g\otimes \g,\wedge^2\g$, 
respectively, hence each of them canonically contains a copy of $\g$. 
Let $I$ be the ideal in $R$ generated by these three copies of $\g$, and 
$A=R/I$. 

The associative algebra $A$ may be interpreted as follows.  
Let $\Pi \g$ denote $\g$ regarded as an odd vector space.
Then $R$ may be thought of as the algebra of regular functions on 
$\Pi\g\times \Pi\g$. 
We have the supercommutator map 
$\lbrace{,\rbrace}: \Pi\g\times \Pi\g\to \g$ 
given by the formula $\lbrace{X,Y\rbrace}=XY+YX$,
where the products are taken in the universal enveloping algebra
(this is a morphism of supermanifolds). The ideal $I$ in $R$ 
is then defined by the relations $\lbrace{X,X\rbrace}=0$, 
$\lbrace{X,Y\rbrace}=0$, $\lbrace{Y,Y\rbrace}=0$, so $A$ 
is the algebra of functions on the superscheme defined by these equations. 

In [CDSW,W], the following conjecture is proposed, and proved
for classical $\g$:

\begin{conjecture}\label{witten}
(i) The algebra $A^\g$ of $\g$-invariants in the algebra $A$ 
is generated by the unique invariant element $S$ of $A$ of degree (1,1)
(namely, $S=Tr|_V(XY)$, where $V$ is a non-trivial irreducible 
finite dimensional representation 
of $\g$).

(ii) $S^g=0$.

(iii) $S^{g-1}\ne 0$. 
\end{conjecture}

\section{New results and conjectures}
Recall that approach of [EK] was based on the consideration 
of the algebra 
$B$ of functions on the superscheme of $X\in \Pi \g$ such that 
$\lbrace{X,X\rbrace}=0$. This algebra was studied by Kostant,  
Peterson, and others (see [K,K1] and references). In particular, it is known 
that as a $\g$-module, $B$ is a direct sum of $2^{{\rm
rk}(\g)}$ non-isomorphic simple 
$\g$-modules $V_{\frak a}$ parametrized by abelian ideals
${\frak a}$ in a Borel subalgebra ${\frak b}$ of $\g$. Namely, if
${\frak a}$ 
is such an ideal of dimension $d$ then it defines (by taking its top exterior 
power) a nonzero vector $v_{\frak a}$ 
in $\wedge^d\g$ (defined up to scaling). This vector generates an irreducible
submodule in $\wedge^d\g$ (with highest weight vector $v_{\frak a}$), 
and the sum of these submodules is a complement to the kernel 
of the projection $\wedge \g\to B$. In fact, this is the unique 
invariant complement, because in each degree $d$ 
it is the eigenspace of the quadratic Casimir $C$ with eigenvalue $d$ (the
largest possible eigenvalue on $\wedge^d\g$; here the Casimir
is normalized so that $C|_{\g}=1$). Thus the (direct) sum 
of $V_{\frak a}$ is canonically identified with $B$ as a $\g$-module. 

The algebra $A^\g$ is obtained from $B\otimes B$ by 
taking a quotient by the additional relation $\lbrace{X,Y\rbrace}=0$ 
and then taking invariants. In other words, $A=(B\otimes
B)^\g/L$, where 
$L$ is the space of invariants in the ideal $\widetilde L$ in
$B\otimes B$ 
given by the relation $\lbrace{X,Y\rbrace}=0$. 

Thus, let us look more carefully at the algebra $E:=(B\otimes
B)^\g$. From the above 
we see that a basis of $(B\otimes B)^\g$ is 
given by the elements $z_{\frak a}$, the canonical elements
in $V_{\frak a}\otimes V_{\frak a}^*$. The element $z_{\frak a}$ 
sits in bidegree $(d,d)$, where $d$ is the dimension of ${\frak
a}$. Thus $E$ is purely even (commutative). 

\begin{proposition}\label{cox} 
(Kostant) The Poincar\'e series $P_E(t)$ 
of $E$ with respect to either of the two gradings 
equals 
$$
\prod_{j=1}^r(1-t^{d_i-1})^{-1}
$$ 
modulo $t^g$, where $r$ is the rank of $\g$ and 
$d_i$ are the degrees of the generating invariants of 
$S\g$ ($d_1=2$). Furthermore, the coefficient of $t^g$ 
in $\prod_{j=1}^r(1-t^{d_i-1})^{-1}-P_E(t)$ is positive.  
\end{proposition}

\begin{proof} 
The first statement is a restatement 
of the results of [K1] (see e.g. Theorem 
4.20 in [K1]). The second statement can be obtained
from the analysis in [CP]. 
\end{proof} 

Now let us explain how to construct elements of $E$. Let  
$F\in (S\g)^\g[k]$ be an invariant of degree $k$. Using $F$,
we can produce an element $\widehat F$ of $E$ given by $\widehat
F(X,Y)=d^2F(\lbrace{X,Y\rbrace})(X\otimes Y)$. In other words, we view $F$ as
an element of $S^{k-2}\g\otimes S^2\g$ and apply the first component
to the bracket $\lbrace{X,Y\rbrace}$ and the second component to $X\otimes 
Y$. Clearly, $\widehat F$ has degree $k-1$ with respect to $X$ 
and the same degree with respect to $Y$.

\begin{proposition} \label{hat} 
Let $F,H$ be two homogeneous invariants in $S\g$ of positive degrees. 
Then $\widehat{FH}=0$. 
\end{proposition}

\begin{proof} We first show that 
for an invariant $F$ of positive degree, 
$F(\lbrace{X,Y\rbrace})=0$ and 
$dF(\lbrace{X,Y\rbrace})(X)=
dF(\lbrace{X,Y\rbrace})(Y)=0$. 
It is well known that the ring $(S\g)^\g$ is generated by 
functions $F(z)=
Tr_V(z^n)$, where $V$ is runs over irreducible representations of 
$\g$ and $n\ge 2$ is an integer. Thus it suffices to prove the claim for such 
$F$. But for such $F$ it is straightforward to check the claim by an easy 
direct computation, using the relations $X^2=Y^2=0$. 

Now, $d^2(FH)=d^2(F)H+Fd^2H+dF\otimes dH+dH\otimes dF.$ 
This implies the proposition. 
\end{proof} 

Let $p_1,...,p_r$ be homogeneous generators of $(S\g)^\g$
(of degrees $d_1=2,...,d_r$). Let $P$ be the span of $\widehat p_i$;
it follows from Proposition \ref{hat} that it is independent of 
the choice of $p_i$. 

\begin{conjecture}\label{c1} 
$P$ generates $E$. 
\end{conjecture}

Evidence for this conjecture is given by 

\begin{proposition}\label{words} 
Let $w$ be a word in $X$ and $Y$ 
and $Q(X,Y)=Tr_V(w(X,Y))$ be an element of $E$.
Then $Q$ belongs to $P$.   
\end{proposition}

\begin{proof} It is easy to see using the relations $X^2=Y^2=0$ that 
$Q=\widehat F$ for some invariant $F$. The rest follows 
from Proposition \ref{hat}.
\end{proof}

\begin{corollary} 
Conjecture \ref{c1} holds for classical Lie algebras.
\end{corollary}

\begin{proof}
For classical Lie algebras, by Weyl's 
fundamental theorem of invariant theory, any invariant 
of $X$, $Y$ is a polynomial of elements of the form 
$\Tr_V(w(X,Y))$, where $V$ is the vector representation. 
Thus, by Proposition \ref{words}, Conjecture 
\ref{c1} holds. 
\end{proof} 

\begin{proposition}\label{cox2}
If Conjecture \ref{c1} is true then 
the lowest degree of a relation among $\widehat p_i$ is $g$. 
\end{proposition}

\begin{proof}
Since the degrees of $p_i$ are $d_i$,
the proposition follows from Proposition \ref{cox}. 
\end{proof}

\begin{proposition}\label{cox3}
 Conjecture \ref{c1} implies part (i) of Conjecture \ref{witten}.
\end{proposition}

\begin{proof}
The proposition follows from the fact that 
for $i>1$ the element $\widehat p_i$ obviously lies in the ideal $L$ in $E$.
\end{proof}

\begin{conjecture}\label{c2}
There exists a relation of degree $g$ among $\widehat p_i$ 
which contains the term 
$\widehat p_1^g$ with a nonzero coefficient. 
\end{conjecture}

\begin{proposition}\label{cox4}
 Conjecture \ref{c2} implies part (ii) of Conjecture \ref{witten}.
\end{proposition}

\begin{proof} Clear. \end{proof} 

\begin{conjecture}\label{c3}
Conjecture \ref{c1} holds, and
the ideal $L$ is generated by $\widehat p_2,...,\widehat p_r$.  
\end{conjecture}

\begin{proposition}\label{cox5}
 Conjecture \ref{c3} implies part (iii) of Conjecture \ref{witten}.
\end{proposition}

\begin{proof} Clear from Proposition \ref{cox}. \end{proof} 

\begin{proposition}\label{23cl}
Conjectures \ref{c2} and \ref{c3} hold for classical Lie algebras. 
\end{proposition}

\begin{proof}
Because of Proposition \ref{words} and Weyl's fundamental theorem of 
invariant theory, Conjecture \ref{c3} 
holds for classical Lie algebras. 
Indeed, by Weyl's theorem any element of $L$ is of the 
form $\Tr_V(w(X,Y)\lbrace{X,Y\rbrace})$,
which implies the desired statement.
Thus, conjecture 
\ref{c2} for classical Lie algebras 
is equivalent to the statement that $S^g=0$, which 
was proved by an explicit computation in [W]. 
Thus in particular we obtain a completely algebraic proof of 
Conjecture \ref{witten} for $SO(N)$. 
\end{proof}

{\bf Remark.} It follows from 
from [EK] that Conjectures \ref{c1},\ref{c2},\ref{c3} 
also hold for $G_2$. 

{\bf Remark.} Let us give a very short proof of 
Conjecture \ref{c2} for $sl(n)$ (cf. also [W]). Let 
$f_n$ be the polynomial such that for $y_i:=x_1^i+...+x_n^i$ 
one has $y_{n+1}=f_n(y_1,...,y_n)$. 
We start with an equation satisfied for any (even)
matrix $Z$: 
$$
\Tr(Z^{n+1})=f_n(\Tr(Z),...,\Tr(Z^n)), 
$$
It is easy to check that 
the coefficient of $y_1^{n+1}$ in $f_n$ is nonzero (it is
$(-1)^n/n!$). Hence the coefficient 
of $y_1^{n-1}y_2$ is also nonzero
(to cancel terms of the form $x_1^2x_2...x_n$). Let us set 
$Z=XY+\xi X+\eta Y$, where $\xi,\eta$ are additional odd
variables, and look at the $\xi\eta$ terms on both sides. 
Then the term $\Tr(Z)^{n-1}\Tr(Z^2)$ (and no other terms) 
yields $\Tr(XY)^n$. This implies the desired statement.  

{\bf Acknowledgements.} The author thanks V. Kac, D. Kazhdan, 
B. Kostant, and A. Postnikov for useful discussions. 
 
\centerline{\bf References}
\vskip .05in

[CDSW]
F. Cachazo, M. R. Douglas, N. Seiberg, E. Witten,
Chiral Rings and Anomalies in Supersymmetric Gauge Theory,
hep-th/0211170.

[CP] 
P. Cellini and P. Papi, 
ad-nilpotent ideals of a Borel subalgebra III,
math.RT/0303065.

[EK] 
P. Etingof and V. Kac, 
On the Cachazo-Douglas-Seiberg-Witten conjecture for simple Lie algebras,
math.QA/0305175.

[K]
B. Kostant,  
On $\bigwedge{\frak g}$ for a semisimple Lie algebra ${\frak g}$, as
an equivariant 
module over the symmetric algebra $S({\frak g})$.
Analysis on homogeneous spaces and representation theory of Lie
groups, 
Okayama--Kyoto (1997), 129--144, Adv. Stud. Pure Math., 26, 
Math. Soc. Japan, Tokyo, 2000.

[K1] B. Kostant, 
Powers of the Euler product and commutative subalgebras of a complex 
simple Lie algebra, math.GR/0309232.  

[W] E. Witten, 
Chiral Ring Of Sp(N) And SO(N) Supersymmetric Gauge Theory In Four Dimensions
hep-th/0302194.

\end{document}